\def\Bbb#1{{\bf #1}}

\def\fnote#1{\footnote}
\def\blacksquare{\hbox{\vrule width 4pt height 4pt depth 0pt}}

\def\cwleftpar#1#2{\leftskip #1 \rightskip #2 plus 1fill}
\def\cwrightpar#1#2{\leftskip #1 plus 1fill \rightskip #2}
\def\cwcenterpar#1#2{\leftskip #1 plus 1fill \rightskip #2 plus 1fill}
\def\cwfullpar#1#2{\leftskip#1\rightskip#2}

\def\cwoutdent#1#2{\llap{\hbox to #1{#2 \hss}}\ignorespaces}
\def\cwparbegin#1#2#3#4#5{
	\ifcase #1 \cwleftpar{#2}{#3}
	\or \cwrightpar{#2}{#3}
	\or \cwcenterpar{#2}{#3}
	\else \cwfullpar{#2}{#3}\fi
	\ifcase #4 \baselineskip = 1.5\baselineskip
	\or \baselineskip = 2\baselineskip
	\or \baselineskip = 3\baselineskip
	\else \baselineskip = 1\baselineskip\fi
	\ifdim #5 > 0in \else \noindent \fi
	\noindent\ignorespaces}
\documentclass{article}
\begin{document}
%------------------------------------------------------------------------
%                               ChiWriter FOOTER
%------------------------------------------------------------------------
\advance \vsize by -1\baselineskip
\def\makefootline{
{\vskip \baselineskip \noindent \folio                                 \par
}}
%------------------------------------------------------------------------
%                               Your Document
%------------------------------------------------------------------------
\vspace*{2ex}
\noindent {\Huge Parallel Transports in Tensor Spaces\\[0.4ex]
Generated by Derivations of\\[0.4ex] Tensor Algebras}

\vspace*{2ex}

\noindent Bozhidar Zakhariev Iliev
\fnote{0}{\noindent $^{\hbox{}}$Permanent address:
Laboratory of Mathematical Modeling in Physics,
Institute for Nuclear Research and \mbox{Nuclear} Energy,
Bulgarian Academy of Sciences,
Boul.\ Tzarigradsko chauss\'ee~72, 1784 Sofia, Bulgaria\\
\indent E-mail address: bozho@inrne.bas.bg\\
\indent URL: http://theo.inrne.bas.bg/$\sim$bozho/}

\vspace*{2ex}

{\bf \noindent Published: Communication JINR, E5-93-1, Dubna, 1993}\\[1ex]
\hphantom{\bf Published: }
http://www.arXiv.org e-Print archive No.~math.DG/0502009\\[2ex]

2003 PACS numbers: 02.40.Ma, 02.40.Vh, 04.90.+e\\[2ex]

\noindent
{\small
The \LaTeXe\ source file of this paper was produced by converting a
ChiWriter 3.16 source file into
ChiWriter 4.0 file and then converting the latter file into a
\LaTeX\ 2.09 source file, which was manually edited for correcting numerous
errors and for improving the appearance of the text.  As a result of this
procedure, some errors in the text may exist.
}\\[2ex]

	\begin{abstract}
The (parallel linear) transports in tensor spaces generated by derivations
of the tensor algebra along paths are axiomatically described. Certain their
properties are investigated. Transports along paths defined by derivations of
the tensor algebra over a differentiable manifold are considered.
	\end{abstract}\vspace{3ex}

 {\bf 1. INTRODUCTION}

\medskip
In [6] we have axiomatically described flat linear transports over a given differentiable manifold and it turns out that the set of these transports coincides with the one consisting of parallel transports generated by flat linear connection. In the present investigation, on this basis, the main ideas of [6] are appropriately generalized and we come to a class of linear transports (along paths in tensor bundles over a manifold) which contains not only the parallel transports generated by arbitrary linear connections, but also the ones generated by arbitrary derivations (along paths).

In Sect. 2 the transports generated by derivations along paths, called here $S$-transports, are axiomatically defined on the basis of restriction of flat linear transports along paths; also some basic properties of the $S$-transports are presented. Sect. 3 is devoted to different relationships between the (axiomatically defined$) S$-transports and derivations along paths. The main result here is the equivalence between the sets of these objects in a sense that to any $S$-transport there can be put into correspondence a unique derivation along paths and vice versa. In Sect. 4 the problem for the connections between $S$-transports and derivations of the tensor algebra over a differentiable manifold is investigated. In particular, it is shown how to any derivation the can be put into correspondence an $S$-transport which in certain "nice" cases turns out to be unique. Also it is shown how some widely used concrete transports along paths (curves), such
as parallel, Fermi-Walker, Fermi etc., can be obtained as special cases from the general $S$-transport. The paper ends with some concluding remarks in Sect. 5.

\medskip
\medskip
 {\bf 2. RESTRICTION ON THE PATH OF FLAT LINEAR\\
 TRANSPORT AND ITS AXIOMATIC DESCRIPTION.\\
	S-TRANSPORTS}

\medskip
 In this work $M$ denotes a real, of class $C^{1}$differentiable manifold
[8]. The tensor space of type $(p,q)$ over $M$ at $x\in M$ will be written as
$T^{p,q}_{x}(M)$. By definition $T_{x}(M):=T^{1,0}_{x}$
and $T^{*}_{x}(M) :=T^{0,1}_{x}(M)$ are, respectively, the tangent and
cotangent spaces to $M$ and $T^{0,0}_{x}(M):={\Bbb R}$. By $J$ and $\gamma :J
\to M$ we denote, respectively, an arbitrary real interval and a path in M.

Let $M$ be endowed with a flat linear transport $L [6]$ and the path $\gamma
:J \to M$ be without self-intersections. Along $\gamma $ we define the map
$S^{\gamma }$ by
\[
 S^{\gamma }:(s,t)\mapsto S^{\gamma }_{s \to t}:=L_{\gamma (s) \to \gamma
(t)},\quad (s,t)\in J \times J.
 \qquad (2.1)
\]
 From the basic properties of
flat linear transports (see definition 2.1 from [6]), we easily derive:
 \[
 S^{\gamma }_{s \to t}(T^{p,q}_{\gamma (s)}(M))\subseteq T^{p,q}_{\gamma
(t)}(M),\quad  s,t\in J,\qquad (2.2)
 \]
\[
  S^{\gamma }_{s \to t}(\lambda A+\mu
A^\prime )=\lambda S^{\gamma }_{s \to t}A+\mu S^{\gamma }_{s \to t}A^\prime ,
\ \lambda ,\mu \in {\Bbb R},\  A,A^\prime \in T^{p,q}_{\gamma (s)}(M),
  \qquad  (2.3)
 \]
 \[
 S^{\gamma }_{s \to t}(A_{1}\otimes A_{2})
 =(S^{\gamma }_{s \to t}A_{1})\otimes
(S^{\gamma }_{s \to t}A_{2}), A_{a}\in T^{p_{a}}_{\gamma (s)}(M),\ a=1,2,
 \qquad (2.4)
  \]
 \[
   S^{\gamma }_{s \to t}\circ C=C\circ S^{\gamma }_{s \to t},\qquad
(2.5)
 \]
 \[
   S^{\gamma }_{t \to r}\circ S^{\gamma }_{s \to t}=S^{\gamma }_{s \to
r},\quad  r,s,t\in J,\qquad (2.6)
 \]
 \[
  S^{\gamma }_{s \to s}=id,
  \qquad (2.7)
 \]
  where $C$
 is any contraction operator and id is the  identity  map (in this case of
 the tensor algebra at $\gamma (s))$.

If for some $s,t\in J$, we have $\gamma (s)=\gamma (t) ($if $s\neq t$ this
 means that $\gamma $ has selfintersection(s)), then from $eq. (2.6)$ of ref.
 [6], we get
\[
  S^{\gamma }_{s \to t}=id\quad if\ \gamma (s)=\gamma (t),
\quad  s,t\in J,\qquad
 (2.8)
 \]
 which characterizes the flat case considered in [6] (see also the
 remark after the proof of proposition 2.1).

{\bf Definition 2.1.} An $S$-transport (along paths) in $M$ is a map $S$
which to any path $\gamma :J \to M$ puts into correspondence a map $S^{\gamma
}, S$-transport along $\gamma $, such that $S^{\gamma }:(s,t)\mapsto
S^{\gamma }_{s \to t}$for $(s,t)\in J\times J$, where the map $S^{\gamma }_{s
\to t}$, an $S$-transport along $\gamma $ from $s$ to $t$, maps the tensor
algebra at $\gamma (s)$ into the tensor algebra at $\gamma (t)$ and satisfies
$(2.2)-(2.7)$.

Above we saw that to any flat linear transport over $M$ there corresponds an
 $S$-transport. The opposite is not, generally, true in a sense that if
 $\gamma :J \to M$ joints $x$ and $y, x,y\in M$, i.e. $x=\gamma (s)$ and
 $y=\gamma (t)$ for some $s,t\in J$, then the map
\[
  L_{x \to y}=S^{\gamma }_{s \to t}, \gamma (s)=x \gamma (t)=y,\qquad
 (2.9)
 \]
  generally, depends on $\gamma $ and is not a flat linear transport
 from $x$ to y.

{\bf Proposition 2.1.} The map (2.9) is a flat linear transport from $x$ to $y$ iff in (2.9) the $S$-transport along every path $\gamma :J \to M$ from $s$ to $t$ depends only on the initial and final points $\gamma (s)$ and $\gamma (t)$, respectively, but not on the path $\gamma $ itself.

{\it Proof.} Let for some $S$-transport the map (2.9) be a flat
linear transport from $x$ to y. As the flat linear transport is a parallel transport generated by a flat linear connection [6], it does not depend on the path along which it is performed [10]. Hence $S^{\gamma }_{s \to t}$depends only on $\gamma (s)=x$ and $\gamma (t)=y$ but not on $\gamma $.

And vice versa, if $S^{\gamma }_{s \to t}$depends only on $\gamma (s)$ and $\gamma (t)$ but not on $\gamma $, then from $(2.3)-(2.7)$ it follows that the map (2.9) satisfies $(2.1)-(2.6)$ from [6], and hence (2.9) is a flat linear transport from $x$ to y.\blacksquare

{\bf Remark.} If the $S$-transport along a "product of paths" is a composition of the $S$-transports along the corresponding constituent paths, then it can be proved that (2.9) defines a flat linear transport iff the $S$-transport in it satisfies (2.8).

A lot of results concerning flat linear transports over $M$ have corresponding analogs in the theory of $S$-transports. Roughly speaking, this transferring of results may be done by replacing points in $M$ with numbers in $J$ if the former do not denote arguments of tensor fields primary defined on $M$; in the last case the points form $M$ must be replaced by the corresponding points from $\gamma (J)\subset $M. In particular, this is true for propositions $2.1-2.5$ from [6], the analogs of which for $S$-transports will be presented below as propositions $2.2-2.6$, respectively. The corresponding proofs will be omitted as they can, evidently, be obtained mutatis mutandis from the ones given in [6].

Let $^{p}_{q}S$ be the restriction of some $S$-transport on a tensor bundle of type $(p,q)$, i.e. $^{p}_{q}S:\gamma \mapsto ^{p}_{q}S^{\gamma },{ } ^{p}_{q}S^{\gamma }:(s,t)\mapsto ^{p}_{q}S^{\gamma }_{s \to t}:= :=S^{\gamma }_{s \to t}$ $_{T^{p,q}_{\gamma (s)}}$. Evidently $^{p}_{q}S$ satisfies $(2.2)-(2.7) ($with $^{p}_{q}S$ instead of $S)$, so the set of the maps $\{^{p}_{q}S, p,q\ge 0\}$ is equivalent to the map S. Hence $S$ splits into the $S$-transports $^{p}_{q}S$ acting independently in different tensor bundles.

{\bf Proposition 2.2.} The linear map $^{p}_{q}S^{\gamma }_{s \to
t}:T^{p,q}_{\gamma (s)}(M) \to T^{p,q}_{\gamma (t)}(M), s,t\in J$ satisfies
(2.6) and (2.7) if and only if there exists linear isomorphisms
$^{p}_{q}F^{\gamma }_{s}:T^{p,q}_{\gamma (s)}(M) \to V, s\in J, V$ being a
vector space, such that
\[
 ^{p}_{q}S^{\gamma }_{s \to t}={\bigl(}^{p}_{q}F^{\gamma }_{t}
{\bigr)}^{-1}{\bigl(}^{p}_{q}F^{\gamma }_{s}{\bigr)}, \quad s,t\in J.
 \qquad (2.10)
\]

 {\bf Proposition 2.3.} If for $^{p}_{q}S^{\gamma }_{s \to
t}$the representation (2.10) holds (see the previous proposition) and $\mathbf{
V}$ is isomorphic with $V$ vector space, then the representation
\[
 ^{p}_{q}S^{\gamma }_{s \to t}={\bigl(}^{p}_{q}\mathbf{ F}^{\gamma
}_{t}{\bigr)}^{-1}{\bigl(}^{p}_{q}\mathbf{ F}^{\gamma }_{s}{\bigr)},
\qquad s,t\in J,\qquad (2.11)
\]
 where $^{p}_{q}\mathbf{ F}^{\gamma }_{s}:T^{p,q}_{\gamma (s)}(M)
\to \mathbf{ V}, s\in J$ are isomorphisms, is true iff there exists isomorphism
$D^{\gamma }:V \to \mathbf{ V}$ such that
\[
^{p}_{q}\mathbf{ F}^{\gamma }_{s}= D^{\gamma }\circ ^{p}_{q}F^{\gamma }_{s},
 \qquad s\in J.  \qquad (2.12)
\]
 Hence on any
fixed tensor bundle every $S$-transport along $\gamma $ from $s$ to $t$
decomposes into a composition of maps depending separately on $s$ and $t$ in
conformity with (2.10). The arbitrariness in this decomposition is described
by proposition 2.3 (see $eq. (2.12))$.

Propositions 2.2 and 2.3 are consequences of $(2.2), (2.6)$ and (2.7). Now we shall also take into account and $(2.3)-(2.5)$.

Putting in $(2.4) A_{1}=A_{2}=1\in {\Bbb R}$, we get $S^{\gamma }_{s \to
t}1=1$ which, due to (2.3), is equivalent to
\[
S^{\gamma }_{s \to t}\lambda =\lambda ,
\quad \lambda \in {\Bbb R}.\qquad(2.13)
\]
 Let $\{E_{i}\mid _{\gamma (s)}\}$ and $\{E^{i}\mid _{\gamma (s)}\}$
be dual bases, respectively, in $T_{\gamma (s)}(M)$ and $T^{*}_{\gamma
(s)}$along $\gamma :J \to M, s\in $J. (The Latin indices run from 1 to
$n:=\dim(M)$ and henceforth we assume the summation rule from 1 to $n$ over
repeated indices.)

From (2.2) it follows that there exist uniquely defined functions
$H^{i}_{.j}(t,s;\gamma )$ and $H^{.j}_{\hbox{i.}}(t,s;\gamma ), s,t\in J$
such that
\[
S^{\gamma }_{s \to t}(E_{j}\mid _{\gamma (s)})=H^{i}_{.j}(t,s;\gamma
)E_{i}\mid _{\gamma (t)},\qquad (2.14a)
\]
\[
 S^{\gamma }_{s \to t}(E^{j}\mid _{\gamma (s)})=H^{.j}_{\hbox{i.}}(t,s;\gamma
)E^{i}\mid _{\gamma (t)}.\qquad (2.14b)
\]
By means of $(2.4), (2.5)$ and (2.13)
it can easily be shown (cf. the derivation of $eq. (2.13)$ from [6]) that
\[
H^{i}_{.k}(t,s;\gamma )H^{.k}_{\hbox{j.}}(t,s;\gamma )=\delta
^{i}_{j}\qquad (2.15)
\]
 which, in a matrix notation reads
\[
H^{i}_{.k}(t,s;\gamma )  H^{.k}_{\hbox{j.}}(t,s;\gamma ) ={\Bbb I}= \delta
^{i}_{j} ,\qquad (2.15^\prime )
\]
 where $\delta ^{i}_{j}$are the Kronecker
deltas and as a first matrix index the superscript is considered.

From (2.14) and (2.3) we conclude that $H^{i}_{.j}(t,s;\gamma )$ and $H^{.j}_{\hbox{i.}}(t,s;\gamma )$ are components of bivectors (two-point vectors) inverse to one another [10], respectively, from $T_{\gamma (t)}(M)\otimes T^{*}_{\gamma (s)}(M)$ and $T^{*}_{\gamma (t)}(M)\otimes T_{\gamma (s)}(M)$. The following proposition shows that they uniquely define the action of the $S$-transport on any tensor.

{\bf Proposition 2.4.} If
$T=T^{i_{1}\ldots i_p}_{j_{1}\ldots j_q} E_{i_{1}}\mid _{\gamma
(s)}\otimes \cdot \cdot \cdot \otimes E_{i_{p}}\mid _{\gamma (s)}\otimes
\otimes E^{j_{1}}\mid _{\gamma (s)}\otimes \cdot \cdot \cdot \otimes
E^{j_{q}}\mid _{\gamma (s)}$, then
\[
 S^{\gamma }_{s \to t}(T)
=
\Bigl( \prod_{a=1}^{p} H^{k_{a}}_{..i_{a}}(t,s;\gamma ) \Bigr)
\Bigl( \prod_{b=1}^{q} H^{..j_{b}}_{l_{b}}(t,s;\gamma ) \Bigr)
T^{i_{1}\ldots i_p}_{j_{1}\ldots j_q}
\]
\[
E_{k_{1}}\mid _{\gamma (t)}\otimes \cdot \cdot \cdot \otimes E_{k_{p}}\mid
_{\gamma (t)}\otimes E^{l_{1}}\mid _{\gamma (t)}\otimes \cdot \cdot \cdot
\otimes E^{l_{q}}\mid _{\gamma (t)}\qquad (2.16)
\]

If $\{E_{i}\mid _{x}\}$ and
$\{e_{i}\}$ are bases in $T_{x}(M)$ and $V$ respectively, then the matrix
elements of $^{1}_{0}F^{\gamma }_{s}, s\in J$ are defined by
$^{1}_{0}F^{\gamma }_{s}(E_{j}\mid _{\gamma (s)})=F^{i}_{.j}(s;\gamma
)e_{i}$. So, if we put $F(s;\gamma ):= F^{i}_{.j}(s;\gamma ) $, then from
(2.10) for $p=q+1=1$ and (2.16), we get
\[
 H(t,s;\gamma ):= H^{i}_{.k}(t,s;\gamma ) =F^{-1}(t;\gamma )F(s;\gamma ),
\quad s,t\in J. \qquad (2.17)
\]

This matrix will be called the matrix of the
considered $S$-transport.

 Evidently (see proposition 2.3), in (2.17) the matrix $F(s;\gamma )$ is
defined up to a constant along $\gamma $ left multiplier, i.e. up to a change
\[
 F(s;\gamma ) \to D^{\gamma }\cdot F(s;\gamma ), \det(D^{\gamma })\neq
0,\infty .\qquad (2.18)
\]

 {\bf Proposition} ${\bf 2}{\bf .}{\bf 5}{\bf .} A$
map $S^{\gamma }_{s \to t}$of the tensor algebra at $\gamma (s)$ into the
tensor algebra at $\gamma (t)$ is an $S$-transport from $s$ to $t$ along
$\gamma $ iff in any local basis its action is given by (2.16) in which the
matrices $ H^{i}_{.k}(t,s;\gamma ) $ and $ H^{.k}_{\hbox{j.}}(t,s;\gamma ) $
are inverse to one another, i.e. (2.15) holds, and the decomposition (2.17)
is valid.

{\bf Proposition 2.6.} Every differentiable manifold admits $S$-transports.

\medskip
\medskip
 {\bf 3. THE EQUIVALENCE BETWEEN S-TRANSPORTS\\
			AND DERIVATIONS ALONG PATHS}
\nopagebreak

\medskip
Let in $M$ an $S$-transport $S$ be given along paths, $\gamma :J \to M$ be
$a C^{1}$path and $T$ be $a C^{1}$tensor field on $\gamma (J)$. To $S$ we
associate a map ${\cal D}$ such that ${\cal D}:\gamma \mapsto {\cal
D}^{\gamma }$, where ${\cal D}^{\gamma }$maps the $C^{1}$tensor fields on
$\gamma (J)$ into the tensor fields on $\gamma (J)$ and
\[
({\cal D}^{\gamma }T)(\gamma (s))
:={\cal D}^{\gamma }_{s}T
:=\lim_{\epsilon \to 0} \frac{1}{\varepsilon}
\bigl[  {\bigl(}
S^{\gamma }_{s+\epsilon  \to s}T(\gamma (s+\epsilon ))-
T(\gamma (s)){\bigr)}\bigr] \qquad (3.1)
\]
 for $s,s+\epsilon \in $J.

Henceforth for the limit in (3.1) to exist, we suppose $S$ to
be of class $C^{1}$in a sense that such is its matrix $H(t,s;\gamma )$ with
respect to $t ($or, equivalently, to $s$; see below (3.4) and (3.5)).

  As a consequence of $(2.7) eq. (3.1)$ can be written also as
\[
{\cal D}^{\gamma }_{s}T
= \Bigl\{
\frac{\partial }{\partial \varepsilon}
 [S^{\gamma }_{s+\epsilon  \to s}
{\bigl(}T(\gamma (s+\epsilon )){\bigr)}]
\Bigr\}\Big|_{\epsilon =0}\qquad (3.1^\prime )
\]

{\bf Proposition 3.1.} The
map ${\cal D}^{\gamma }$is a derivation of the restriction of the tensor
algebra over $M$ on $\gamma (J)$.

{\it Proof.} It can easily be verified that from $(3.1^\prime ), (2.3),
(2.2), (2.5)$ and (2.4) it follows respectively that ${\cal D}^{\gamma }$is
an ${\Bbb R}$-linear, type preserving map of the restriction of the tensor
algebra over $M$ on $\gamma (J)$ into itself which commutes with the
contractions and obeys the Leibnitz rule, i.e. on $\gamma (J)$, we have
\[
{\cal D}^{\gamma }(\lambda A+\mu A^\prime )=\lambda {\cal D}^{\gamma }A +
\mu {\cal D}^{\gamma }A^\prime , \lambda ,\mu \in {\Bbb R},\qquad (3.2a)
\]
\[
{\cal D}^{\gamma }\circ C=C\circ {\cal D}^{\gamma },\qquad (3.2b)
\]
\[
{\cal D}^{\gamma }(A\otimes B)=({\cal D}^{\gamma }A)\otimes B+A\otimes
({\cal D}^{\gamma }B),\qquad (3.2c)
\]
where $C$ is a contraction operator, $A, A^\prime $ and $B$ are $C^{1}$tensor
fields on $\gamma (J)$, the types of A and $A^\prime $ being the same. By
definition [8] this means that ${\cal D}^{\gamma }$is a derivation of the
mentioned restricted algebra.\blacksquare

{\bf Proposition 3.2.} For every path $\gamma :J \to M, s,t\in J$ and
$S$-tran-sport $S^{\gamma }$along $\gamma $ there is valid the identity
\[
{\cal D}^{\gamma }_{t}\circ S^{\gamma }_{s \to t}\equiv 0,\qquad (3.3)
\]
${\cal D}^{\gamma }_{t}$ being defined by $S^{\gamma }$ through (3.1).

 {\it Proof.} The identity (3.3) follows from $(3.1^\prime )$ and
(2.6).\blacksquare

{\bf Proposition 3.3.} If $T$ is $a C^{1}$of type $(p,q)$ tensor field on
$\gamma (J)$ with local components $T^{i_{1}}_{j_{1}}$in a basis $\{E_{i}\}$
defined on $\gamma (J)$, then the local components of ${\cal D}^{\gamma }T$
at $\gamma (s), s\in J$, i.e. of ${\cal D}^{\gamma }_{s}T$, are
\[
{\bigl(}{\cal D}^{\gamma }_{s}T{\bigr)}
^{i_{1}\ldots i_p}_{j_{1}\ldots j_q}
={\bigl(}
\frac{d}{ds} T^{i_{1}\ldots i_p}_{j_{1}\ldots j_q}
{\bigr)}\big|_{\gamma (s)}
\]
\[
+ \sum_{a=1}^{p}\Gamma ^{i_{a}}_{..k}(s;\gamma )
T^{i_{1}\ldots i_{a-1}ki_{a+1}\ldots i_p}_{j_{1}\ldots ]j_q}(\gamma(s))
\]
\[
 -\sum_{b=1}^{a}\Gamma ^{k}_{.j_{b}}(s;\gamma )
T^{i_{1}\ldots i_p}_{j_{1}\ldots j_{b-1}kj_{b+1}\ldots j_q}(\gamma
(s)),\qquad (3.4)
\]
where
\[
\Gamma ^{i}_{.j}(s;\gamma )
:= \frac{\partial H^{i}_{.j}(s,t;\gamma) }{\partial t} \Big|_{t=s}
=-  \frac{\partial H^{i}_{.j}(t,s;\gamma) }{\partial t} \Big|_{t=s}.
\qquad (3.5)
\]

{\it Proof.} This result is a direct consequence of $(3.1^\prime ),
(2.16)$ and (2.15).\blacksquare

{\bf Remark 1.} If we define $\Gamma ^{\gamma }(s):= \Gamma
^{i}_{.j}(s;\gamma ) $, then, due to (2.17), the second equality in (3.5) is
a corollary of
\[
\Gamma ^{\gamma }(s)
= \frac{\partial H(s,t,\gamma}{\partial t}\Big|_{t=s}
= F^{-1}(s;\gamma ) \frac{\partial F(s;\gamma}{\partial s} .
\qquad (3.5^\prime )
\]

 {\bf Remark 2.} If the vector fields $E_{i}$are
defined on $\gamma (J)$ such that $E_{i}\mid _{\gamma (s)}, s\in J$ is a
basis at $\gamma (s)$, then from (3.4) it follows that
\[
 {\cal D}^{\gamma }E_{i}=(\Gamma ^{\gamma })^{k}_{.i}E_{k}=\Gamma
^{k}_{.i}(s;\gamma )E_{k}\qquad (3.5^{\prime\prime})
\]
and, vice versa, if we  define $\Gamma ^{\gamma }$by the expansion
$(3.5^{\prime\prime})$, then, from proposition 3.1 (see (3.2)) and (3.1), we
easily get (3.4).

If $\nabla $ is a covariant differentiation with local components $\Gamma
^{i}_{.jk}$and $\gamma $ is a $C^{1}$path with a tangent vector field
$\dot\gamma$, then the comparison of the explicit form of
$(\nabla_{\dot\gamma}T)(\gamma (s)) (cf.  [8, 10])$ with (3.4) shows that the
derivation ${\cal D}^{\gamma }$is a generalization of the covariant
differentiation $\nabla _{\cdot }$along $\gamma $. Evidently, ${\cal
D}^{\gamma }$reduces to a covariant differentiation along $\gamma $ iff
\[
\Gamma ^{i}_{.j}(s;\gamma )=\Gamma ^{i}_{.jk}(\gamma (s))\dot\gamma^{k}(s),
\qquad (3.6)
\]
 where $\Gamma ^{i}_{.jk}$ are components of some
covariant differentiation $\nabla $.

{\bf Lemma 3.1.} The change $\{E_{i}\mid _{\gamma (s)}\} \to \{E_{i^\prime
}\mid _{\gamma (s)}=A^{i}_{i^\prime }(s)E_{i}\mid _{\gamma (s)}\}, s\in J$
leads to the transformation of $\Gamma ^{i}_{.j}(s;\gamma )$ into
\[
\Gamma ^{i^\prime }_{..j^\prime }(s;\gamma )=A^{i^\prime
}_{i}(s)A^{j}_{j^\prime }(s)\Gamma ^{i}_{.j}(s;\gamma )+A^{i^\prime
}_{i}(s)(dA^{i}_{j^\prime }(s)/ds),\qquad (3.7)
\]
 where $[ A^{i^\prime }_{i}] :=[ A^{j}_{j^\prime }(s)] ^{-1}$.

{\it Proof.} $Eq. (3.7)$ is a simple corollary of (3.5) and the fact that
$H^{i}_{.j}(t,s;\gamma )$ are components of a tensor from $T_{\gamma
(t)}(M)\otimes T^{*}_{\gamma (s)}($see (2.14)).\blacksquare

{\bf Proposition 3.4.} If in any basis $\{E_{i}\mid _{\gamma (s)}\}, s\in J$
along $\gamma :J \to M$ there are given functions $\Gamma ^{i}_{.j}(s;\gamma
)$ which, when the basis is changed, transform in conformity with (3.8), then
there exists a unique $S$-transport along $\gamma $ which generates $\Gamma
^{\gamma }(s):= \Gamma ^{i}_{.j}(s;\gamma ) $ through $(3.5^\prime )$ and the
matrix $H(t,s;\gamma )$ of which is
\[
 H(t,s;\gamma )=Y(t,s_{0};-\Gamma ^{\gamma })[Y(s,s_{0}-\Gamma ^{\gamma
})]^{-1},\quad  s,t\in J. \qquad (3.8)
\]
Here $s_{0}\in J$ is fixed and $Y=Y(s,s_{0};Z), Z$ being a continuous matrix function of $s$, is the unique
solution of the initial-value problem
\[
 dY/ds=ZY,\quad Y\mid _{s=s_{0}}={\Bbb I}.\qquad (3.9)
\]

 {\bf Remark.} The existence and uniqueness of the solution of (3.9) can be
found, for instance, in [4].

{\it Proof.} At first we shall prove that for a fixed $\Gamma ^{\gamma }eq.
(3.8)$ gives the unique solution of $(3.5^\prime )$ with respect to
$H(t,s;\gamma )=F^{-1}(t;\gamma )F(s;\gamma )$. In fact, using
$dF^{-1}/ds=-F^{-1}(dF/ds)F^{-1}$, we see $(3.5^\prime )$ to be equivalent to
$dF^{-1}(s;\gamma )/ds=-\Gamma ^{\gamma }(s)F^{-1}(s;\gamma )$, the general
solution of which with respect to $F^{-1}$, due to (3.9) is $F^{-1}(s;\gamma
)=Y(s,s_{0};-\Gamma ^{\gamma })D(\gamma )$, where $s_{0}\in J$ is fixed and
$D(\gamma )$ is a nondegenerate matrix function of $\gamma $. Substituting
the last expression of $F^{-1}$into (2.17), we get (3.8) (which does not
depend either on $D(\gamma )$ or on $s_{0}; cf$. proposition 2.3 and the
properties of $Y ($see [4])).

It can easily be shown that, as a consequence of the transformation law (3.7), the elements of $H(t,s;\gamma )$ are components of a tensor from $T_{\gamma (t)}(M)\otimes T^{*}_{\gamma (s)}$. Hence, by proposition 2.4 (see also (2.15)), they define an $S$-transport, the local action of which is given by (2.16). Due to the above construction of $H(t,s;\gamma )$ this $S$-transport is the only one that generates $\Gamma ^{\gamma }$, as given by $(3.5^\prime ).\blacksquare $

So, the definition of an $S$-transport along $\gamma :J \to M$ is equivalent to the definition of functions $\Gamma ^{i}_{.j}(s;\gamma ), s\in J$, in a basis $\{E_{i}\}$ along $\gamma $, which have the transformation law (3.7). For this reason we shall call the functions defined by (3.5) components of the $S$-transport.

Proposition 3.1 shows that to any $S$-transport there corresponds, according to $(3.1), a$ derivation of the restriction of the tensor algebra over $M$ on any curve $\gamma (J)$, i.e. along the path $\gamma :J \to $M. The next proposition states that any such derivation can be obtained in this way.

{\bf Proposition 3.5.} The map ${\cal D}^{\gamma }$of the restriction of the tensor algebra over $M$ on $\gamma (J)$ into itself is a derivation iff there exists an $S$-transport along $\gamma :J \to M ($which is unique and) which generates ${\cal D}^{\gamma }$by means of $eq. (3.1)$.

{\it Proof.} If ${\cal D}^{\gamma }$is defined by some $S$-transport through (3.1), then by proposition 3.1 it is a derivation. And vice versa, let ${\cal D}^{\gamma }$be a derivation, i.e. to be a type preserving and satisfying
(3.2). If we define $\Gamma ^{\gamma }:= \Gamma ^{i}_{.j}(s;\gamma ) $ by $(3.5^{\prime\prime})$, it is easily verified that the transformation law (3.7) holds for $\Gamma ^{i}_{.j}(s;\gamma )$ and, consequently (see proposition 3.4), there exists a unique $S$-transport for which $\Gamma ^{i}_{.j}(s;\gamma )$ are local components in the used basis. The derivation corresponding, in conformity with (3.1), to this $S$-transport has an explicit action given by the right hand side of (3.4) and hence it coincides with ${\cal D}^{\gamma }$as the latter has the same explicit action (as a consequence of $(3.5^{\prime\prime})).\blacksquare $

So, any derivation of the restricted along a path tensor algebra is generated by a unique $S$-transport along this path through (3.1). The opposite statement is almost evident and it is expressed by

{\bf Proposition 3.6.} For any $S$-transport along a path $\gamma :J \to M$ there exists a unique derivation ${\cal D}^{\gamma }$of the restricted along $\gamma $ tensor algebra which generates its components $\Gamma ^{\gamma }(s):= \Gamma ^{i}_{.j}(s;\gamma ) $ in a basis $\{E_{i}\mid _{\gamma (s)}\}, s\in J$ through equation $(3.5^{\prime\prime})$.

{\it Proof.} The existence of ${\cal D}^{\gamma }$for a given $S$-transport is evident: by proposition 3.1 the map ${\cal D}^{\gamma }$, defined by (3.1), is a derivation along $\gamma $ and by proposition 3.3 it generates the components of the $S$-transport by (3.5), or, equivalently, by $(3.5^{\prime\prime})$. If $^\prime {\cal D}^{\gamma }$is a derivation along $\gamma $ with the same property, then from (3.2) and (3.4), we get ${\cal D}^{\gamma }T=^\prime {\cal D}^{\gamma }T$ for every $C^{1}$tensor field $T$, i.e. $^\prime {\cal D}^{\gamma }={\cal D}^{\gamma }.\blacksquare $

Thus we see that there is a one-to-one correspondence between $S$-transports along paths and derivations along paths, i.e. maps ${\cal D}$ such that for every $\gamma :J \to M$, we have ${\cal D}:\gamma \mapsto {\cal D}^{\gamma }$, where ${\cal D}^{\gamma }$is a derivation of the restricted on $\gamma (J)$ tensor algebra over M.

\medskip
\medskip
 {\bf 4. LINEAR TRANSPORTS ALONG PATHS DEFINED BY\\
 	DERIVATIONS. SPECIAL CASES}

\medskip
Let ${\cal D}$ be a derivation along paths (see the end of Sect. 3). The
identity (3.3) gives the following way for defining the $S$-transport
corresponding to ${\cal D}$, in conformity with proposition 3.5. If $\gamma
:J \to M$, then we define the needed $S$-transport along $\gamma $ by the
initial-value problem
\[
{\cal D}^{\gamma }_{t}\circ S^{\gamma }_{s \to t}=0,
\quad  S^{\gamma }_{s \to t}\mid _{t=s}=id,
 \qquad (4.1)
\]
where for any tensor field $T {\cal D}^{\gamma }_{t}(T):=({\cal
D}^{\gamma }T)(\gamma (t))$. In fact, the initial-value problem (4.1) has a
unique solution with respect to $S^{\gamma }_{s \to t}[4]$ and it is easily
verified that it is an $S$-transport along $\gamma $ from $s$ to $t$, whose
matrix is given by (3.8) in which $\Gamma ^{\gamma }$is defined by
$(3.5^{\prime\prime})$. By proposition 3.3 this means that the derivation
corresponding to this $S$-transport (see (3.1) and proposition 3.1) coincides
with ${\cal D}^{\gamma }$, i.e. it generates ${\cal D}^{\gamma }$.

Hence (4.1) naturally generates the $S$-transport corresponding to some derivation along paths ${\cal D}$. The opposite is true up to a left multiplication with nonzero functions. i.$e$ it can be proved that if the $S$-transport is fixed, then all derivations along paths for which (4.1) is valid have, according to (3.3), the form $f\cdot {\cal D}$, in which $f$ is a nonvanishing scalar function, and ${\cal D}$ is the derivation generated by (3.1) from the given $S$-transport.

If $D$ is a derivation of the tensor algebra over $M$, i.e. if it is a type preserving map satisfying (3.2) on $M$, then along a path $\gamma :J \to M$ to it there may be assigned an $S$-transport along $\gamma $ in the following way.

According to proposition 3.3 from $ch. I$ in [8], there exist a
unique vector field $X$ and a unique tensor field $S$ of type (1,1) such that
\[
 D=L_{X}+S,\qquad (4.2)
\]
where $L_{X}$is the Lie derivative along $X$ and
$S$ is considered as  a derivation of the tensor algebra over $M [8]$. On the
opposite, for every vector field $X$ and tensor field $S_{X}$of type (1,1),
which may depend on $X$, the equation
\[
D=L_{X}+S_{X},\qquad (4.2^\prime )
\]
 defines a derivation of the tensor algebra over M.

Let a derivation $D$ with the decomposition $(4.2^\prime )$ be given and $a
C^{1}$path $\gamma :J \to M$ be fixed. Let in some neighborhood of $\gamma
(J) a$ vector field $V$ be defined such that on $\gamma (J)$ it reduces to
the vector field tangent to $\gamma $, i.e. $V\mid _{\gamma (s)}=$ $(s), s\in
$J. We define ${\cal D}^{\gamma }$to be the restriction of $D\mid _{X=V}$on
$\gamma (J)$, i.e. for every tensor field $T$, we put
\[
 ({\cal D}^{\gamma }T)(\gamma (s)):={\cal D}^{\gamma }_{s}T:=\{(D\mid
_{X=V})(T)\}(\gamma (s)).\qquad (4.3)
\]
 It is easily seen that the map ${\cal
D}^{\gamma }$thus defined is a derivation along $\gamma $. Of course,
generally, ${\cal D}^{\gamma }$depends on the values of $V$ outside the set
$\gamma (J)$. However the most interesting and geometrically valuable case is
that when ${\cal D}^{\gamma }$depends only on the path $\gamma $ but not on
the values of $V$ at points not lying on $\gamma (J)$. In this case the map
${\cal D}:\gamma \mapsto {\cal D}^{\gamma }$is a derivation along the paths
in $M$, therefore according to the above scheme, it generates an
$S$-transport along the same paths and, consequently, in this way the
derivation $D$ generates an $S$-transports. This $S$-transport depends only
on D. Vice versa, if ${\cal D}^{\gamma }$depends on the values of $V$ outside
$\gamma (J)$, then to construct a derivation along paths ${\cal D}:\gamma
\mapsto {\cal D}^{\gamma }$we have to fix in a neighborhood of every path
$\gamma  a$ vector field $V^{\gamma }$such that $V^{\gamma }\mid _{\gamma
(s)}=$ $(s), s\in J$ and to put $({\cal D}^{\gamma }T)(\gamma (s)):=(D$
$_{X=V^{\gamma }}T)(\gamma (s))$. Analogously, using the above described
method, from ${\cal D}$ we can construct an $S$-transport, but ${\cal D}$, as
well as the $S$-transport, will depend generally not only on $D$ but on the
family of vector fields $\{V^{\gamma }: \gamma :J \to M\}$ that is to a great
extent arbitrary.

Hence to any derivation there corresponds at least one $S$-transport. Without going into details of the problem when this transport depends or not on the family $\{V^{\gamma }\}$, we shall present some examples important from the practical point of view examples in which the transport does not depend on that family. (All known to the author and used in the mathematical and physical literature transports, based on affine connections over a manifold, are of this kind, i.$e$ they do not depend on $\{V^{\gamma }\}.)$

Let $M$ be an $L_{n}$space, i.e. it is to be endowed with a linear
connection $\nabla $. As for any vector field $X \nabla _{X}$is a derivation,
in conformity with $(4.2^\prime )$, it admits the representation
\[
 \nabla _{X}=L_{X}+\Sigma (X),\qquad (4.4)
\]
 where $\Sigma (X)$ is a tensor
field of type (1,1). It can easily be proved that the local components of
$\Sigma (X)$ in $\{E_{i}\}$ are $(\Sigma (X))^{i}_{.j}= =(\nabla
_{E_{i}}X)^{i}-T^{i}_{.jk}X^{k}$in which $T^{i}_{.jk}$are the torsion tensor
components.

Expressing $L_{X}$from (4.4) and substituting it into $(4.2^\prime )$, we
get the unique decomposition of any derivation in the form
\[
 D=\nabla _{X}+S_{X}-\Sigma (X)\qquad (4.5)
\]
 which turns out to be useful
for comparison of the general $S$-transport with concrete linear transports,
based on linear connections, of tensors along paths. For instance, using it,
the method presented above for generating $S$-transports from derivation, and
the definitions of the concrete transports mentioned below, given in the
cited there references, one can easily prove the following proposition (cf.
proposition 4.1 from [7]).

{\bf Proposition 4.1.} The $S$-transport generated by a derivation $D$, with
a decomposition (4.5), is reduced:

a) in a space $L_{n}$with a linear connection to the parallel transport $[8,
10]$ for $S_{X}=\Sigma (X)$;

$b)$ in an Einstein-Cartan space $U_{n}$to the Fermi-Walker transport $[5,
11]$ for $S_{X}=\Sigma (X)-2L$;

$c)$ in a Riemannian space $V_{n}$to the Fermi transport [11] for
$S_{X}=\Sigma (X)-2\mathbf{ L}$;

$d)$ in a Riemannian space $V_{n}$to the Truesdell transport [12] for
$S_{X}=\theta \cdot \delta $;

$e)$ in a Riemannian space $V_{n}$to the Jaumann transport [9] for
$S_{X}=\Sigma (X)-\omega $,

where $\delta $ is the unit tensor (with Kronecker symbols as components$),
\theta :=$ $(\nabla _{E_{i}}X)^{i}$is the expansion of $X [5, 9, 12]$, and
$L, \mathbf{ L}$ and $\omega $ are tensor fields of type (1,1) with,
respectively, the following covariant components:
$L_{ij}:=-h^{k}_{i}h^{l}_{j}X_{[k;l]}+(h_{il}T^{l}_{.jk}X^{k})_{[ij]}+V_{[i;j]}$, in which $h_{ij}:=g_{ij}-X_{i}X_{j}/(g_{kl}X^{k}X^{l}), g_{ij}$being the metric components, $X^{k}:=(\nabla _{E_{l}}X)^{k}$, and $(\ldots  )_{[ij]}$means antisymmetrization (e.g., $X_{[i;j]}:=(X_{i;j}-X_{j;i})/2); \mathbf{ L}_{ij}:=X_{i}X_{j;k}X^{k}$for $g_{ij}X^{i}X^{j}=-1$; and $\omega _{ij}:=(X_{i;j}-X_{i;k}X^{k})_{[ij]}$.

The list of concrete transports in proposition 4.1 can be extended to include the $M$-transport [2], the Lie transport $[5, 10]$, the modified Fermi-Walker and the Frenet-Serret transports [1] etc., but this is only a technical problem which does not
change the main idea that by an appropriate choice of $S_{X}$and the the application of the above-described procedure one can obtain a number of useful transports of tensors along paths.

\medskip
\medskip
 {\bf 5. COMMENTS}

\medskip
In the present paper we have considered the axiomatic approach to transports of tensors along paths the generated by derivations along paths, called here $S$-transports. This is done on the basis of axiomatic description of the parallel transport generated by flat linear connections and gives possibilities for further generalizations to be a subject of other papers. On this ground a number of properties of the $S$-transports are derived. We have proved that in a natural way to any $S$-transport there corresponds a unique derivation along paths and vice versa. If one considers general derivations of the tensor algebra, then this correspondence still exists but it is, generally, not unique and needs in this sense an additional investigation.

It can be shown that the $S$-transports are special cases of the parallel transports in fibre bundles $[3, 8]$. Elsewhere we shall see that the theory developed here can be generalized so as to include and these parallel transports as its special case.

\medskip
\medskip
 {\bf ACKNOWLEDGEMENTS}

\medskip
The author expresses his gratitude to Prof. Stancho Dimiev and Prof. Vl. Aleksandrov (Institute of Mathematics of Bulgarian
Academy of Sciences) for constant interest in this work and stimulating discussions. He thanks Prof. N. A. Chernikov (Joint Institute for Nuclear Research, Dubna, Russia) for the interest in the problems posed in this work.

This research was partially supported by the Foundation for Scientific Research of Bulgaria under contract Grant No. $F 103$.

\medskip
\medskip
 {\bf REFERENCES}

\medskip
1. Dandoloff R., W. J. Zakrzewski, Parallel transport along space curve and related phases, J. Phys., {\bf A}: Math. Gen., vol. $22, 1989, L461-L466$.\par
2. Dixon W. G., Proc. Roy. Soc. Lond., Ser. {\bf A}, vol. ${\bf 2}{\bf 1}{\bf 4}, 1970, 499-527$.\par
3. Dubrovin B., S. P. Novikov, A. Fomenko, Modern Geometry, I. Methods and Applications, Springer Verlag.\par
4. Hartman Ph., Ordinary Differential Equations (John Wiley \& Sons, New York-London-Sydney, $1964), ch$.IV, \S1.\par
5. Hawking S., G. F. R. Ellis, The large Scale Structure of Spacetime, Cambridge Univ. Press, Cambridge, 1973.\par
6. Iliev B. Z., Flat linear connections in terms of flat linear transports in tensor bundles, Communication JINR, $E5-92-544$, Dubna, 1992.\par
7. Iliev B. Z., Linear transports of tensors along curves: General $S$-transport, Comp. Rend. Acad., Bulg. Sci., vol. {\bf 40}, No.$7, 1987, 47-50$.\par
8. Kobayashi S., K. Nomizu, Foundations of Differential Geometry, Vol. 1, Interscience Publishers, New York-London, 1963.
\par
9. Radhakrishna L., L. N. Katkar, T.H. Date, Gen. Rel. Grav., vol. {\bf 13}, No. $10, 1981, 939-946$.\par
10. Schouten J. A., Ricci-Calculus: An Introduction to Tensor Analysis and its Geometrical Applications, $2-nd ed$., Springer Verlag, Berlin-G\"ottingen-Heidelberg, 1954.\par
11. Synge J. L., Relativity: the general theory, North-Holland Publ. Co., Amsterdam, 1960.\par
12. Walwadkar B. B., Gen. Rel. Grav., vol. {\bf 15}, No. $12, 1983, 1107-1114$; Walwadkar B. B., K. V. Vikar, Gen. Rel. Grav., vol. {\bf 16}, No. $1, 1984, 1-7$.

\newpage

\vspace{5ex}
\noindent
 Iliev B. Z.\\[5ex]

\noindent
 Parallel  transports in tensor spaces \\
 generated by derivations of tensor algebras\\[5ex]

The (parallel linear) transports in tensor spaces generated by derivations
of the tensor algebra along paths are axiomatically described. Certain their
properties are investigated. Transports along paths defined by derivations of
the tensor algebra over a differentiable manifold are considered.\\[5ex]

The investigation has been performed at the Laboratory of Theoretical
Physics, JINR.

\end{document}